\documentclass[10pt]{article}

\pdfoutput=1

\usepackage{amsmath,amsthm,amssymb,amscd,fancybox,ifthen,float}
\floatstyle{plain}

\usepackage{graphicx}
\usepackage[all]{xy}
\usepackage[margin=1in]{geometry}
\usepackage{setspace}
\usepackage{verbatim}
\usepackage{caption}
\usepackage{subcaption}
\usepackage{mathpazo}
\usepackage{comment}

\newcommand\bc{\boldsymbol c}
\newcommand\bx{\boldsymbol x}
\newcommand\by{\boldsymbol y}
\newcommand\bn{\boldsymbol n}

\theoremstyle{definition}

\theoremstyle{remark}
\newtheorem{remark}{Remark}

\newcommand\blfootnote[1]{%
  \begingroup
  \renewcommand\thefootnote{}\footnote{#1}%
  \addtocounter{footnote}{-1}%
  \endgroup
}

\usepackage{sectsty}
\sectionfont{\large}
\subsectionfont{\normalsize}

\begin{document}

\onehalfspacing
\numberwithin{equation}{section}

\title{On the efficient representation of the half-space
impedance Green's function for the Helmholtz equation}

\author{Michael O'Neil\blfootnote{Research supported in part by the National Science
Foundation under grant DMS06-02235, the U.S. Department of Energy
under contract DE-FG02-88ER-25053, and the Air Force Office of 
Scientific Research under NSSEFF Program Award FA9550-10-1-0180.},
 Leslie Greengard, and Andras Pataki\\
\small{Courant Institute of Mathematical Sciences, New York
University, New York, NY 10012}}

\numberwithin{equation}{section}
\maketitle

\begin{abstract}
A classical problem in acoustic (and electromagnetic) scattering
concerns the evaluation of the Green's function for the Helmholtz
equation subject to impedance boundary conditions on a half-space.
The two principal approaches used for representing this Green's
function are the Sommerfeld integral and the (closely related)
method of complex images. The former is extremely efficient when
the source is at some distance from the half-space boundary, but
involves an unwieldy range of integration as the source gets closer
and closer. Complex image-based methods, on the other hand, can be
quite efficient when the source is close to the boundary, but they
do not easily permit the use of the superposition principle since
the selection of complex image locations depends on both the source
and the target. We have developed a new, hybrid representation
which uses a finite number of real images (dependent only on the
source location) coupled with a rapidly converging Sommerfeld-like
integral. While our method applies in both two and three
dimensions, we restrict the detailed analysis and numerical
experiments here to the two-dimensional case.
\end{abstract}

\section{Introduction}

\begin{figure}[t]
\centering
\includegraphics[width=5.25in]{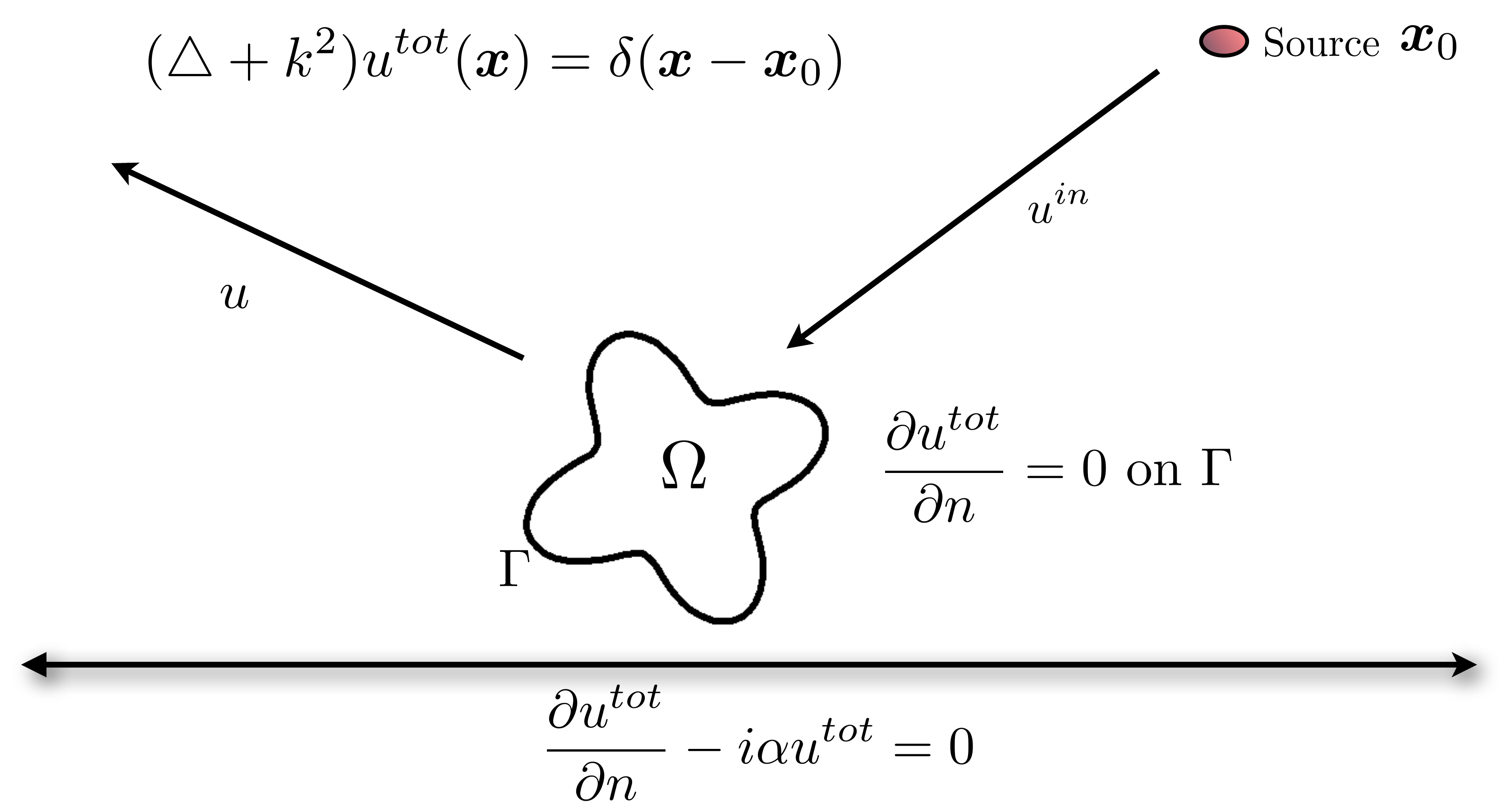}
\caption{Scattering from a \emph{sound-hard} obstacle above an impedance plane.}
\label{fig_scat}
\end{figure}

A number of problems in acoustics (and electromagnetics) involve 
the solution of the Helmholtz equation,  
\begin{equation}\label{eq_helmpde}
(\triangle+ k^2) u(\bx) = f(\bx) \, ,
\end{equation}
in the half-space $P= \{ (x,y) \in {\mathbb R}^2: y > 0 \}$ or 
$S= \{ (x,y,z) \in {\mathbb R}^3: z > 0 \}$, subject to
suitable boundary and radiation conditions.  In acoustics, the
Helmholtz coefficient $k$ is given by $k = \frac{\omega}{c}$, where
$\omega$ is the governing angular frequency (assuming a time-harmonic
motion dependency of $e^{-i\omega t}$) and $c$ is the sound speed.  In
the present paper, we assume $k \in \mathbb C$ is constant throughout
the region of interest, with ${\rm Re}(k) \geq 0$ and ${\rm Im}(k)
\geq 0$.  For concreteness, we concentrate initially on the two-dimensional 
problem of computing the scattered field due to a unit-strength
point source located at $\bx_0=(x_0,y_0)$ in the presence of a
"sound-hard" obstacle over an infinite half-space subject to impedance
boundary conditions (Figure~\ref{fig_scat}).

We let the total field be defined as $u^{tot} = u^{in} + u$, where
$u^{in}$ denotes the (known) incoming field due to the point source
and $u$ denotes the scattered field.  On a sound-hard obstacle
$\Omega$ with boundary $\Gamma$, the total field must satisfy
homogeneous Neumann boundary conditions.  Since the scattered field
involves no sources outside $\Omega$, it must satisfy the homogeneous
Helmholtz equation
\begin{equation}\label{eq_helmpdehom}
  (\triangle+ k^2) u(\bx) = 0 \
\end{equation}
for $ \bx \in P \setminus \Omega$ .  On the obstacle boundary
$\Gamma$, we have
\begin{equation}\label{eq_neumannbc}
  \frac{\partial u}{\partial n} = - \frac{\partial u^{in}}{\partial n},
\end{equation}
where $\frac{\partial }{\partial n}$ is the outward normal derivative.
Finally, on the interface, we assume a standard impedance
condition on the total field of the form:
\begin{equation}\label{eq_imped}
\frac{\partial u^{tot}}{\partial n} -  i \alpha u^{tot} = 0.
\end{equation}
Since the interface is the $x$-axis, we have $\frac{\partial
}{\partial n} = -\frac{\partial }{\partial y}$. In
physically-motivated problems, an impedance condition is typically
used to approximate a more complicated wave/surface interaction,
such as scattering from a rough surface, an underlying porous
medium, a complicated surface coating, etc. (see
\cite{attenborough, cw_horoshenkov}). In many applications,
$\alpha=\beta k$, with $0 \leq \beta \leq 1$, in which case any
dissipation is due entirely to the imaginary part of $k$. The
parameter $\beta$ in this context is called the surface admittance.
In other cases, the physical model introduces dissipation of some
other kind, resulting in a complex valued $\alpha$, even when $k$
is real. For the purposes of this paper, we will assume that
$\alpha \in \mathbb C$, with ${\rm Re}(\alpha) \geq 0$, ${\rm
Im}(\alpha) \geq 0$, $|\alpha| \leq |k|$, and leave aside any
further discussion of the modeling. The Green's function analysis
of the present paper can be generalized to other values of
$\alpha$, but we restrict our attention to $\alpha$ in the
indicated range for the sake of simplicity. A second simplification
is that we only consider the case of constant $\alpha$ (i.e. we do
not permit $\alpha$ to vary along the length of the half-space
interface). There is a substantial literature on impedance problems
and we mention only a few relevant papers which also discuss the
computation of the corresponding Green's function. These include
\cite{caiyu,chandler-wilde,colton_kress,cw_hothersall,
nedelec-impedance,ochmann,sarabandi,taraldsen,vanderpol}.

Returning now to the scattering problem 
(\ref{eq_helmpdehom}, \ref{eq_neumannbc}, \ref{eq_imped}),
an ansatz for the solution is to
represent the total field as 
\begin{equation} \label{eq-slprep}
  u^{tot}(\bx) = \int_\Gamma g_{k,\alpha} (\bx,\by) \, \sigma(\by) \,
  d{\by(s)} + u^{in}(\bx), 
\end{equation}
where $s$ is arclength along $\Gamma$, 
$g_{k,\alpha}(\bx,\bx_0)$ is the Green's function for the half-space
$P$ with homogeneous impedance boundary conditions, and
$u^{in}(\bx) = g_{k,\alpha}(\bx,\bx_0)$.
Imposing the Neumann
conditions (\ref{eq_neumannbc}) 
on $\Gamma$ yields the Fredholm integral equation of the second kind:
\begin{equation}
\label{eq-neu}
  -\frac{1}{2} \sigma(\bx) + \int_\Gamma \frac{\partial }{\partial
    n_x}g_{k,\alpha}(\bx,\by) \, \sigma(\by) \, d{\by(s)} = -
  \frac{\partial }{\partial n_x}g_{k,\alpha}(\bx,\bx_0)
\end{equation}
for $\bx \in \Gamma$, where the integral is interpreted in the
principal value sense. Equation (\ref{eq-neu}) is invertible except
for a countable sequence of spurious resonances $\{k_j\}$.
Resonance-free, but more complicated, representations are well-known
\cite{colton_kress}, which we will not review here, since we are
primarily interest in the question of how to efficiently evaluate the
impedance Green's function $g_{k,\alpha}$ itself. In our examples, we
will always assume $k \notin \{k_j\}$ and that equation
(\ref{eq-neu}) is solvable.  Note that, by using the impedance
Green's function in the integral representation, the infinite
half-space boundary does $\emph{not}$ need to be discretized.  

Algorithms for the computation of $g_{k,\alpha}$ date back to the
classical work of Sommerfeld, Weyl, and Van der Pol
\cite{sommerfeld,weyl,vanderpol}, who developed both what are now
referred to as the Sommerfeld integral and the method of complex
images.  For more recent treatments of this problem, see
\cite{hochmann,digilbert,dinapoli,ochmann,taraldsen,koh}.

The main contribution of the present work is the observation that a
finite number of {\em real} images can accurately capture the
high-frequency components of the  Sommerfeld integral. This leads,
naturally, to a hybrid representation of the Green's function in
terms of a rapidly converging Sommerfeld-type representation,
augmented with $\mathcal O(\log(1/d))$ real images for each source
point that lies a distance $d$ from the impedance interface.  Our
approach is somewhat related to that of Cai and Yu \cite{caiyu},
which also separates low- and high-frequency contributions, but uses
an asymptotic method for the high-frequency components.

The paper is organized as follows. Section~\ref{sec_spectral} gives
a derivation of the classical spectral representation for the free
space Green's function, due to Sommerfeld. In
Section~\ref{sec_impedance}, we discuss Sommerfeld and Van der
Pol's extension of the spectral representation to the case of
impedance boundary conditions for a half-space.
Section~\ref{sec:images} introduces analytical (closed-form)
expressions for the real and complex image representations. In
Section~\ref{sec:hybrid}, we present our new representation that
combines a finite segment of real images in the lower half-space
with a Sommerfeld integral that is rapidly decaying.
Section~\ref{sec-disc} discusses some the details concerning
discretization and quadrature for both the image segment and the
obstacle boundary $\Gamma$, and Section~\ref{sec_numerical}
contains several numerical experiments which demonstrate the
effectiveness of the scheme. Lastly, in
Section~\ref{sec_conclusions}, we discuss the extension of the
method to the three-dimensional case, to layered media, and to the
Maxwell equations - all areas for future research.

\section{Spectral representation of the Green's function}
\label{sec_spectral}

The solution $g_k$ to the Helmholtz equation 
\begin{equation}
  \label{eq_helm}
  (\triangle+ k^2) g_k(\bx) = \delta(\bx-\bx_0),
\end{equation}
in an infinite homogeneous medium is referred to as the free-space
Green's function, where $\bx = (x,y) \in \mathbb R^d$, $d$ is the
underlying dimension, and $\delta(\bx-\bx_0)$ represents the Dirac
delta function centered at $\bx_0$.  It is well known that
\[
  g_k(\bx,\bx_0) = \left\{ \begin{array}{cc}
  \frac{i}{4} H^{(1)}_0(k|\bx-\bx_0|) & {\rm for}\ d=2, \\
  \quad & \quad \\ 
  -{\frac{e^{ik|\bx-\bx_0|}}{4\pi |\bx-\bx_0|}} & {\rm for}\ d=3,
  \end{array} \right.
\]
where $H^{(1)}_0$ denotes the zeroth-order Hankel function of the
first kind.  These Green's functions satisfy the outgoing Sommerfeld
radiation condition
\begin{equation}
\lim_{r \to \infty} r^{(d-1)/2} \left( \frac{\partial}{\partial r}
g_k(\bx,\bx_0) - ik \, g_k(\bx,\bx_0) \right),
\end{equation}
where $r=|\bx-\bx_0|$.

A continuous spectral representation of the Green's functions can be obtained
by taking the Fourier transform of equation~(\ref{eq_helm}). In two 
dimensions,  the Green's function can then be written as
\begin{equation}
g_k(\bx,\bx_0)=\frac{1}{4\pi^2} \int_{-\infty}^\infty 
\int_{-\infty}^\infty
\frac{e^{i(\lambda_x(x-x_0)+\lambda_y(y-y_0))}}
{\lambda_x^2 + \lambda_y^2 -k^2} \,
d\lambda_x  \, d\lambda_y.
\end{equation}
Evaluating the
integral in $\lambda_y$ via contour deformation
yields the expansion in plane waves (often called the Sommerfeld integral): 
\begin{equation}
\label{eq_gkhelm3}
g_k(\bx,\bx_0) =
\frac{1}{4\pi} \int_{-\infty}^\infty \frac{e^{-\sqrt{\lambda^2-k^2}|y-y_0|}}
{\sqrt{\lambda^2-k^2}} \, e^{i\lambda (x-x_0)} \, d\lambda.
\end{equation}
Due to the central role of this formula in scattering theory,
there has been much effort devoted to its numerical evaluation. 
We do not give a comprehensive review of the various schemes available. 
They are largely based on contour deformation into the second and fourth
quadrants in the complex $\lambda$-plane in order to avoid
the square-root singularity in the denominator.
In our numerical calculations, we make use of a hyperbolic tangent contour, 
\[ \lambda(t) = t - i \tanh(t) 
\]
for $t \in (-\infty,\infty)$, as in 
\cite{barnett_greengard} and the trapezoidal rule on the interval 
$-t_{max} \leq t \leq t_{max}$ for some $t_{max}$. 
Assuming the integrand has vanished
(to high precision) at the endpoints $\pm t_{max}$, this results in 
a spectrally accurate quadrature scheme.
The difficulty in computing the Sommerfeld integral is clear from 
(\ref{eq_gkhelm3}); when $|y-y_0|$ is small, the integrand 
is slowly decaying and the range of integration prohibitively
large. 

\section{The impedance problem}
\label{sec_impedance}

In the context of the half-space problem,
we need an analytic representation
of the response to the free-space Green's function that 
enforces the homogeneous impedance condition.
This can be done in either the frequency domain, as in 
(\ref{eq_gkhelm3}), or by introducing an infinite ray of images emanating from the
reflection of the source point across the $x$-axis.
Our method is based on combining these two ideas.

For the spectral approach
\cite{hochmann,digilbert,dinapoli,ochmann,taraldsen,koh,sommerfeld,weyl,vanderpol}, 
we begin by using (\ref{eq_gkhelm3}) to represent
the field induced on the impedance boundary by a single free-space
point source located at $\bx_0 = (x_0,y_0)$.
One can then match Fourier modes
$e^{i\lambda x}$ between the known incoming field and the unknown scattered
field to enforce the desired impedance condition.
More precisely, in the two-dimensional
case, the incoming field for points $\bx = (x,y)$ with $y< y_0$
can be written as
\begin{equation}
  u^{in}_{free}(\bx)=g_k(\bx,\bx_0)=\frac{1}{4\pi} \int_{-\infty}^\infty 
  \frac{e^{\sqrt{\lambda^2-k^2}(y-y_0)}}
  {\sqrt{\lambda^2-k^2}} \, e^{i\lambda(x-x_0)} \, d\lambda.
\end{equation}
Suppose now that we assume the analogous
spectral representation for the scattered field
$u$: 
\begin{equation}
  \label{eq_sighat}
  u(\bx)=\frac{1}{4\pi} \int_{-\infty}^\infty 
  \frac{e^{-\sqrt{\lambda^2-k^2}y}}{\sqrt{\lambda^2-k^2}} \,
  e^{i\lambda x} \, \hat\sigma(\lambda) \, d\lambda,
\end{equation}
where $\hat\sigma$ is an unknown density.  (The
formula~(\ref{eq_sighat}) is the spectral representation of a single
layer potential due to a charge density $\sigma$ on the interface
$y=0$.)  Imposing the impedance boundary condition (\ref{eq_imped}) on
the interface $y=0$, we obtain a simple scalar equation for
$\hat\sigma(\lambda)$:
\begin{equation}
  -e^{-\sqrt{\lambda^2-k^2} y_0} e^{-i\lambda x_0} +
  \hat\sigma(\lambda) -i\alpha \left( \frac{e^{-\sqrt{\lambda^2-k^2}
      y_0}}{\sqrt{\lambda^2-k^2}} e^{-i\lambda x_0}
  +\frac{\hat\sigma(\lambda)}{\sqrt{\lambda^2-k^2}} \right) = 0.
\end{equation}
Solving for the density $\hat\sigma$ yields
\begin{equation}
  \hat\sigma(\lambda)= e^{-\sqrt{\lambda^2-k^2}y_0} \, e^{-i\lambda x_0}
  \left( \frac{\sqrt{\lambda^2-k^2}+i\alpha}{\sqrt{\lambda^2-k^2}-i\alpha}
  \right) \, ,
\end{equation}
and the scattered field $u$ can be evaluated via the representation
\begin{equation}
\label{eq_imsom}
u(\bx)=\frac{1}{4\pi} \int_{-\infty}^\infty 
\frac{e^{-\sqrt{\lambda^2-k^2}(y+y_0)}}
{\sqrt{\lambda^2-k^2}} \, e^{i\lambda (x-x_0)}
\left(
\frac{\sqrt{\lambda^2-k^2}+i\alpha}{\sqrt{\lambda^2-k^2}-i\alpha}\right)
\, d\lambda.
\end{equation}
The full impedance Green's function is then given by
$g_{k,\alpha}=u+u^{in}_{free}$. 
Note that since $y$ and $y_0$ are both
positive, (\ref{eq_imsom}) can be used to efficiently evaluate the scattered
field if {\em either} $y$ or $y_0$ are $\mathcal O (1)$, since the integrand
is exponentially decaying for large $\lambda$.  However, if $y+y_0
\sim \mathcal O (h) \ll 1$, then the size of the integration interval
must be chose to be $\mathcal O ({h}^{-1})$, which can be unreasonably
large. This is the case when both the scatterer and the target are near the interface.

A variety of attempts have been made to introduce more efficient schemes in this regime.
Cai and Yu, for example, added an artificial mollifier to the Sommerfeld representation
and expanded the remaining high frequency components asymptotically
\cite{caiyu}. More common, however, is the use of the method of images, both
real and complex, which we turn to next.

\section{The method of images} \label{sec:images}

The use of image charges to impose a given homogeneous boundary
condition is a well-known technique in classical applied mathematics
\cite{kellogg}.
When solving the half-space problem with homogeneous Dirichlet
boundary conditions, for example, the response to a free-space point source located
at $(x_0,y_0)$ is exactly the field generated by a point source of
equal and \emph{opposite} strength located at $(x_0,-y_0)$.
Similarly, for the homogeneous Neumann problem, the response to a
point source located at $(x_0,y_0)$ is exactly the field generated by
a point source of \emph{equal} strength located at $(x_0,-y_0)$.
Unfortunately, in the case of impedance boundary conditions,
no single image source is sufficient.  However, it is possible to
develop an explicit representation of the impedance Green's function
using an infinite ray of images, starting at the reflected 
point $(x_0,-y_0)$ and continuing vertically down (see
\cite{taraldsen} for a historical overview).

\begin{figure}[t]
\centering
\includegraphics[width=6in]{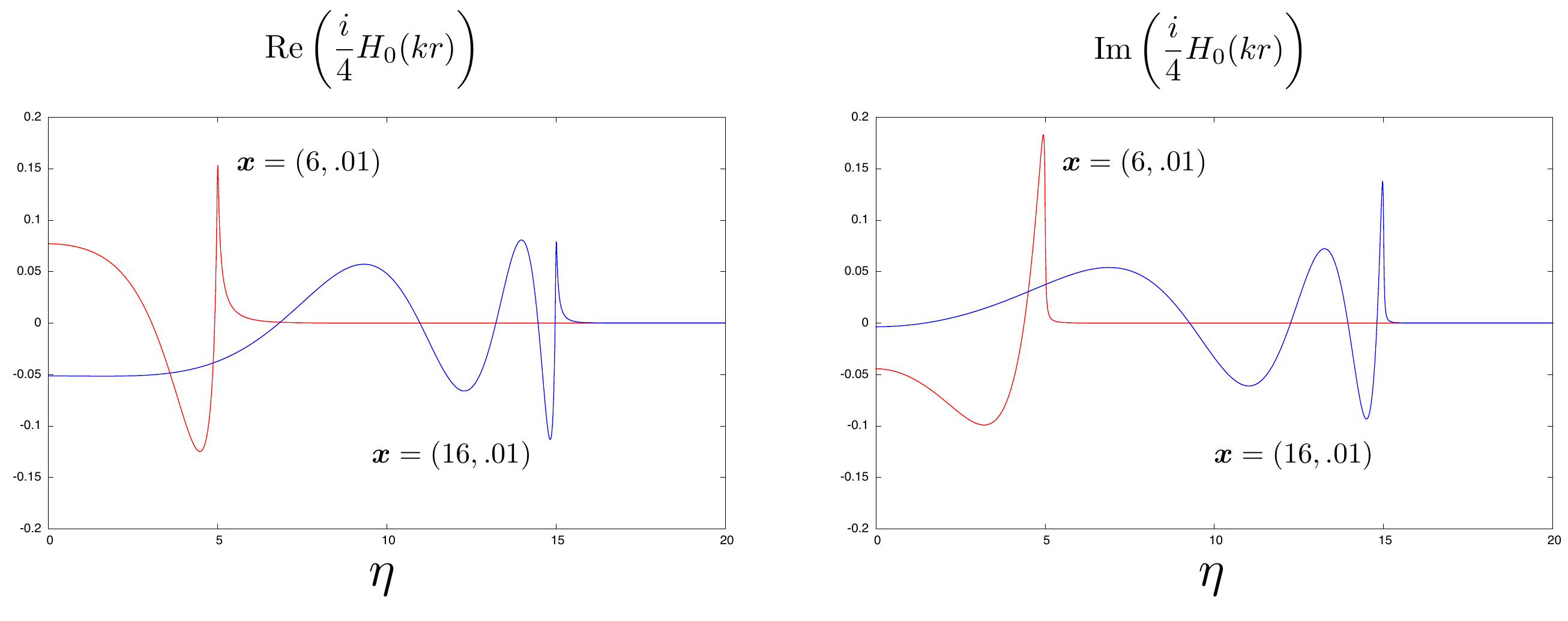}
\caption{For $k=1$ and source location $\bx_0=(1,.01)$, the plot shows
  the behavior of the singularity of $H_0$ and its dependence on the
  lateral location of the source and target. The distance $r$ is given
  by $r~=~\sqrt{(x-x_0)^2+(y+y_0+i\eta)^2}$, corresponding to a
  complex-valued source location.}
\label{fig_complex}
\end{figure}

In this approach, we assume the scattered field $u$ takes the form:
\begin{equation}
\label{eq:image}
 u(\bx) = \int_0^\infty g_k(\bx,\bx_0-(2y_0 +\eta)\hat\by)
\, \tau(\eta) \, d\eta,
\end{equation}
where $\hat\by=(0,1)$ is the unit normal vector in the $y$ direction
and $\tau(\eta)$ is an unknown charge distribution.
Using the Sommerfeld representation for the free-space Green's function $g_k$, 
we may write the scattered field $u$ as
\begin{equation}
\label{eq_linerep}
u(\bx) = \frac{1}{4\pi} \int_0^\infty  \int_{-\infty}^\infty 
\frac{e^{-\sqrt{\lambda^2-k^2}(y+y_0+\eta)}}
{\sqrt{\lambda^2-k^2}} \, e^{i\lambda (x-x_0)} \, 
\tau(\eta)  \, d\lambda \, d\eta.
\end{equation}
In order for $u^{tot}$ to satisfy the impedance boundary condition at
$y=0$, it is straightforward to see that the density $\tau$ must satisfy
the following equation for all $\lambda$:
\begin{equation}
-e^{-\sqrt{\lambda^2-k^2}y_0}
+\int_0^\infty e^{-\sqrt{\lambda^2-k^2}(y_0+\eta)}
\, \tau(\eta) \, d\eta -
i\alpha \left( \frac{e^{-\sqrt{\lambda^2-k^2}y_0}}{\sqrt{\lambda^2 - k^2}}
+ \int_0^\infty \frac{e^{-\sqrt{\lambda^2-k^2}(y_0+\eta)}}
{\sqrt{\lambda^2-k^2} } \, \tau(\eta) \, d\eta 
\right) = 0.
\end{equation}
After some algebra, this reduces to a condition on the Laplace transform
of the image density $\tau$,
\begin{equation}\label{eq_mu}
\begin{split}
\int_0^\infty e^{-\sqrt{\lambda^2-k^2}\eta} \, \tau(\eta) \, d\eta &=
\frac{\sqrt{\lambda^2-k^2} + i\alpha}{\sqrt{\lambda^2-k^2} - i\alpha} \\
&=1+2i\alpha \frac{1}{ \sqrt{\lambda^2-k^2} - i\alpha }.
\end{split}
\end{equation}
This equation can be solved by inspection, using two simple Laplace
transform identities:
\begin{equation}
\begin{split}
\mathcal L [e^{-\beta t}](s) &= \frac{1}{s+\beta}, \\
\mathcal L [\delta(t-t_0)](s) &= e^{-t_0 s},
\end{split}
\end{equation}
where
\[
\mathcal L[f](s) = \int_0^\infty e^{-st} \, f(t) \, dt.
\]
It is easy to see that the
solution to ~(\ref{eq_mu}) is given by
\begin{equation}
\tau(\eta) = \delta(\eta)+2i\alpha e^{i\alpha\eta}.
\end{equation}
The complete image-based formula for the scattered field $u$ is then
\begin{equation}
\label{eq_image}
u(\bx) = g_k(\bx,\bx_0-2 y_0 \hat\by) + 2i\alpha
\int_0^\infty g_k(\bx,\bx_0-(2 y_0 + \eta)\hat\by) \,
 e^{i\alpha\eta} \, d\eta ,
\end{equation}
and once again,
the impedance Green's function can be 
constructed as $g_{k,\alpha}=u+u^{in}_{free}$.
Note that if $\alpha$ is purely real,  the density $\tau(\eta)$ oscillates 
throughout its range,
so that the decay in the integrand comes only from the 
$1/\eta$  decay in the Green's function $g_k$.
To overcome this, the standard solution involves 
complexification of the coordinate $\eta$ (see \cite{taraldsen,vanderpol}).
The simple change of variables $\eta \to i\eta$ in formula
(\ref{eq:image}), for example, yields
\begin{equation}
\label{eq_image2}
u(\bx) = g_k(\bx,\bx_0-2 y_0 \hat\by) - 2\alpha
\int_0^\infty g_k(\bx,\bx_0-(2 y_0 + i\eta)\hat\by) \,
 e^{-\alpha\eta} \, d\eta \, .
\end{equation}
For impedance parameters $\alpha$ with ${\rm Re}(\alpha) > 0$, this
complex image representation enforces exponential decay of the
integrand and avoids the main difficulty in evaluating the integral
(\ref{eq_image}).  Unfortunately, a side-effect of this change of
variables is that the behavior of the integrand is rather complicated
- involving the evaluation of the free-space Green's function over a
range of complex arguments as $\eta$ varies.  Figure~\ref{fig_complex}
illustrates the behavior of the integrand in equation (\ref{eq:image})
at two distinct target locations with the same $y$~value. This
behavior prevents the straightforward design of fast numerical
algorithms that make use of the principal of superposition. In other
words, because of the sensitivity of the singularity in the integrand
to the location of targets and sources, it is difficult to find
robust, efficient, and universal quadratures for the integral in
(\ref{eq_image2}).

\section{A hybrid approach} \label{sec:hybrid}

\begin{figure}[t]
\centering
\includegraphics[width=3in]{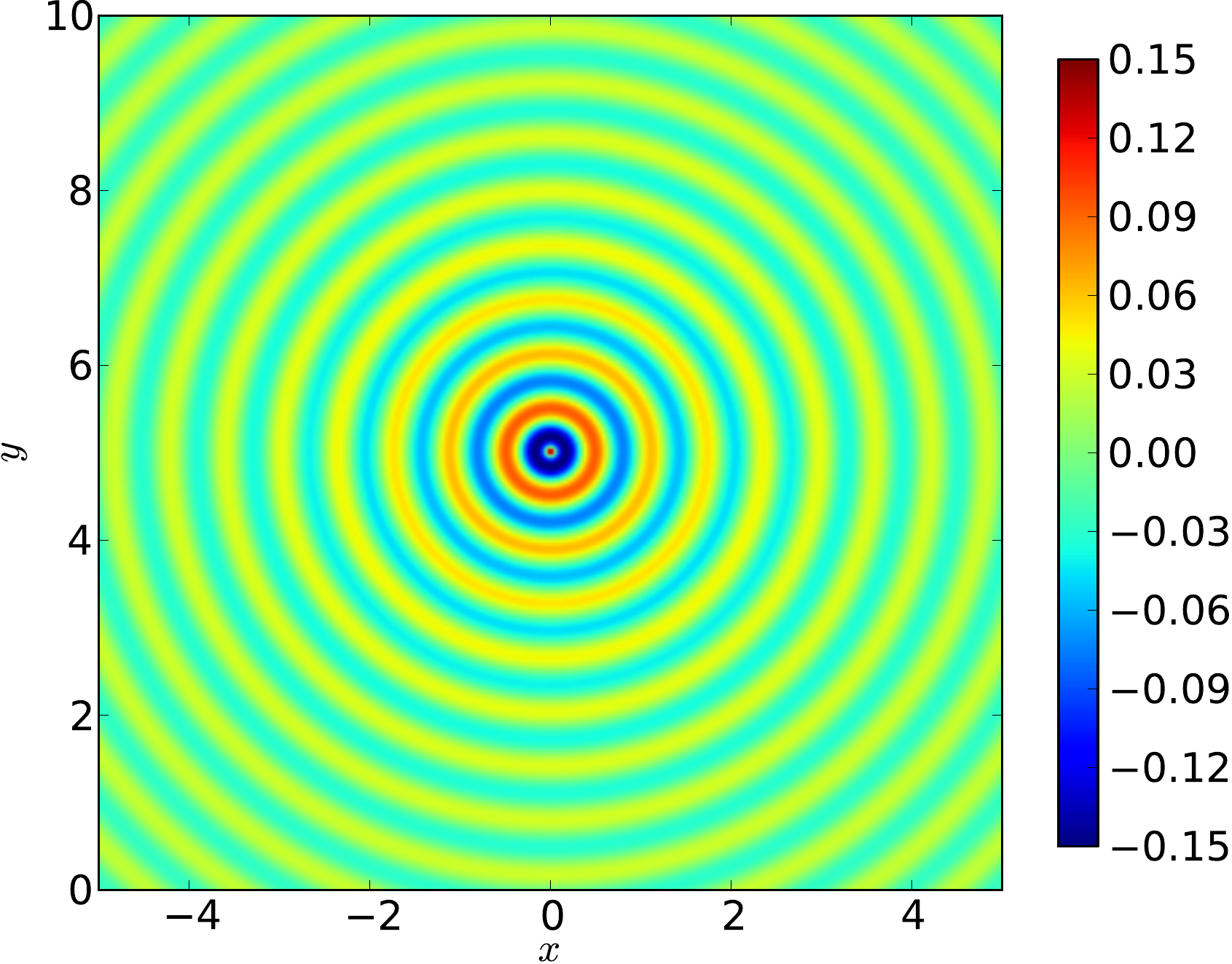}
\quad \includegraphics[width=3in]{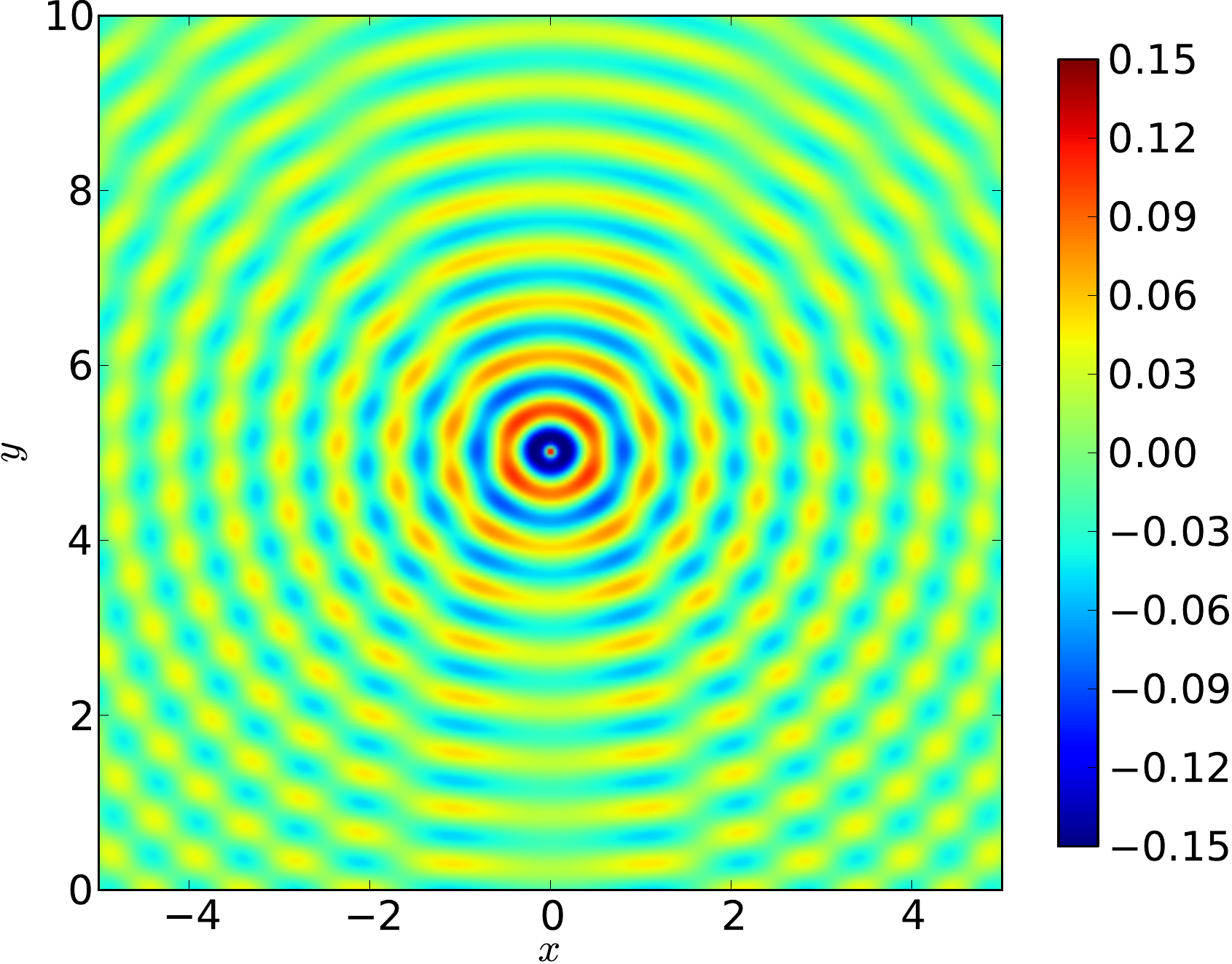}
\caption{For $k=10.2$, the plot on the left shows the
real part of the two dimensional free-space Green's function
and the plot on the right shows the real part of the impedance
Green's function with $\alpha=2.04$. Partial reflection from the
boundary $y=0$ is clearly visible.}
\label{fig_heat}
\end{figure}

It turns out that there is a representation of the impedance Green's 
function which can take advantage of 
both the Sommerfeld integral approach and the method of images.
We begin by reconsidering the {\em real} image formula~(\ref{eq_image}),
and separating the image ray into two parts:
a \emph{near-field} component and a \emph{far-field} component. The
scattered field $u$ from formula~(\ref{eq_image}) is then written as
\begin{equation}
\begin{split}
\label{eq_im2}
u(\bx) =&   \ g_k(\bx,\bx_0-2 y_0 \hat\by) \\
&+2i\alpha
\left( \int_0^{C} g_k(\bx,\bx_0-(2 y_0 -\eta)\hat\by)
 \, e^{i\alpha\eta} \, d\eta +
\int_{C}^\infty \! g_k(\bx,\bx_0-(2 y_0 -\eta)\hat\by)
 \, e^{i\alpha\eta} \, d\eta \right) ,
\end{split}
\end{equation}
where $C$ is a parameter of our choosing.
From equations (\ref{eq_imsom}) and (\ref{eq_im2}), it is straightforward
to show that
\begin{equation}
\label{eq_farfield}
 \int_C^\infty g_k(\bx,\bx_0-2 y_0 \hat\by-\eta\hat\by) \,
 e^{i\alpha\eta} \, d\eta
= \frac{1}{4\pi} \int_{-\infty}^\infty 
\frac{e^{-\sqrt{\lambda^2-k^2}(y+y_0)}} {\sqrt{\lambda^2-k^2} }
\, \frac{e^{-(\sqrt{\lambda^2-k^2}-i\alpha)C}}
{\sqrt{\lambda^2-k^2} - i\alpha} \, e^{i\lambda (x-x_0)} \, d\lambda.
\end{equation}
This is a Sommerfeld-like formula which, for $C \sim \mathcal O (1)$,
decays exponentially fast once $\lambda \gtrsim |k|$, independent of
$y$ and $y_0$.  Therefore, the full impedance Green's function can be
written as
\begin{equation}
\label{eq:hybrid_rep}
\begin{split}
  g_{k,\alpha}(\bx,\bx_0) =& \ g_k(\bx,\bx_0) + 
  g_k(\bx,\bx_0-2 y_0 \hat\by) + 2i\alpha
  \int_0^{C} g_k(\bx,\bx_0-2 y_0 \hat\by-\eta\hat\by) \,
  e^{i\alpha\eta} \, d\eta \\
   &+\int_{-\infty}^\infty 
      \varphi_{k,\alpha,C} (\bx,\bx_0,\lambda) \, d\lambda,
\end{split} 
\end{equation}
where
\begin{equation} \label{varphi_def}
\varphi_{k,\alpha,C} (\bx,\bx_0,\lambda) =
\frac{i\alpha}{2\pi} \,
  \frac{e^{-\sqrt{\lambda^2-k^2}(y+y_0)}} {\sqrt{\lambda^2-k^2} }
  \frac{e^{-(\sqrt{\lambda^2-k^2}-i\alpha)C}}
  {\sqrt{\lambda^2-k^2} -i\alpha} \, e^{i\lambda(x-x_0)} \, .
\end{equation}
Using the previous formula for evaluation, 
Figure~\ref{fig_heat} plots (for comparison) the real part of the free-space
Green's function and the impedance half-space Green's function, for a source
located at $(0,5)$ with $k=10.2$.

\section{Discretization and fast algorithms}
\label{sec-disc}

In our numerical experiments, the impedance Green's function
is discretized and evaluated as follows: 
\begin{equation}\label{eq_sum}
\begin{split}
g_{k,\alpha}(\bx,\bx_0) = &g_k(\bx,\bx_0) + g_k(\bx,\bx_0-2 y_0
\hat\by) + \sum_{j=1}^{Q_I} v_j \, 
g_k(\bx,\bx_0-(2 y_0 + \eta_j) \hat\by) \\ 
&+\sum_{j=1}^{Q_S} w_j \, \varphi_{k,\alpha,C} (\bx,\bx_0,\lambda_j),
\end{split}
\end{equation}
where $\eta_j,v_j$ represents the $j^\text{th}$ quadrature node and
weight used in the evaluation of the first integral
in~(\ref{eq:hybrid_rep}), and $\lambda_j,w_j$ represents the
$j^\text{th}$ quadrature node and weight used in the evaluation of the
second integral in~(\ref{eq:hybrid_rep}).  Here, $Q_I$ and $Q_S$
denote the total number of nodes in the discretizations of the first
and second integrals, respectively.  When the dependence on $\bx_0$ is
needed, we will write $Q_I(\bx_0)$ or $Q_S(\bx_0)$.  The nodes and
weights $\{\eta_j,v_j\}$ are chosen to be $16^\text{th}$ order
Gauss-Legendre nodes on the dyadic subdivisions $[0,2^{-m}C],
[2^{-m}C,2^{-m+1}C],\dots, [2^{-1}C,C]$.  The number of intervals
$(m+1)$ is chosen \emph{a priori} so that the integral is evaluated to
a specified precision $\epsilon$ when $y=0$.  Thus, the number of
nodes $Q_I \approx 16m$.  For a scattering source at $(x_0,y_0)$, the
integral is accurate once the size of the smallest dyadic interval is
$\mathcal O\left(\frac{y_0}{k}\right)$ so that $Q_I(\bx_0) =
\mathcal O\left(\log\frac{k}{y_0}\right)$.  For the Sommerfeld
component, the nodes and weights $\lambda_j,w_j$ are chosen to
correspond to the trapezoidal rule in the real variable $t \in
[-t_{max},t_{max}]$ along the contour
\begin{equation}
    \lambda=t - i \tanh t.
\end{equation}
For $k$ bounded away from zero, this quadrature provides spectral
accuracy due to the decay and smoothness of the integrand in the
Sommerfeld-like contribution, assuming $t_{max}$ is chosen so that
$e^{-(\sqrt{\lambda^2-k^2}-i\alpha)C}$ is small. This is easy to check
and $t_{max} = |k| + 20$ is typically sufficient for 10 digits, with
the total number of nodes required $Q_S$ of the order $\mathcal O(k +
\alpha)$, assuming $(x-x_0) = O(1)$.  Note that the number of nodes
required \emph{does not} increase as the the distance of the source
and target from the interface goes to zero.

\begin{remark}
\label{rem-c}
\onehalfspacing
We have not, as yet, specified $C$ from equation~(\ref{eq_im2}).
In all of our experiments, we simply let $C=1$ independent of the
source location. In fact, when the source is well separated from the
interface (say by a distance of $1/k$, the Sommerfeld integral
converges rapidly enough so that the real images aren't essential.
\end{remark}

\subsection{A basic fast algorithm}
\label{sec-fast}

We now sketch an outline of a fast algorithm for evaluating
the impedance Green's function at $N$ targets due to $M$ sources.
We denote the $M$ point sources~$\{ \bx_m' \}$ and the 
$N$ target locations by $\{ \bx_l \}$. Thus,
\begin{equation}
u(\bx_l) = \sum_{m=1}^M c_m \, g_{k,\alpha}(\bx_l,\bx_m')
\qquad \text{for $l=1,\ldots,N$}.
\end{equation}
We proceed in two steps:
\begin{enumerate}
\item Compute the first three terms on the right
hand side of equation~(\ref{eq_sum}). These are merely sums of free-space
point sources with {\em real} coordinates, 
so their contribution can be computed using a
two-dimensional Helmholtz fast multipole method (FMM) in $\mathcal O(
n \log n ) $ time \cite{gg_software, wideband2d, rokhlin}, where
$n = N + M + M_{image}$. Here, $M_{image} = \sum_{m=1}^M Q_I(\bx_m')$ 
denotes the total number of image points generated for all sources.

\item Compute the contribution of the Sommerfeld-like integral for each
target. A naive calculation would require $\mathcal O(NM)$
operations. However, the representation makes it very easy to exploit
separation of variables.  Denoting by $u_S(\bx_l)$ the Sommerfeld-like
contribution at the point $\bx_l$, and using (\ref{varphi_def}), we
have
\begin{equation}
\begin{split}
u_S(\bx_l) &= 
\sum_{m=1}^M c_m \sum_{j=1}^{Q_S}  \, w_j \,
 \varphi_{k,\alpha,C} (\bx,\bx_m',\lambda_j)  \\
&= \frac{i\alpha}{2\pi}
 \sum_{j=1}^{Q_S} \,  w_j
\frac{e^{-\sqrt{\lambda_j^2-k^2}y_l}} {\sqrt{\lambda_j^2-k^2} }
\, \frac{e^{-(\sqrt{\lambda_j^2-k^2}+i\alpha)C}}
{\sqrt{\lambda_j^2-k^2} +i\alpha} \, e^{-i\lambda_j x_l} \,
\rho(\lambda_j,k),
\end{split}
\end{equation}
where
\[
\rho(\lambda_j,k) = \sum_{m=1}^M c_m \,
e^{-\sqrt{\lambda_j^2-k^2}y_m'} \, e^{i\lambda x_m'}.
\]
The evaluation of $\rho(\lambda_j,k)$ requires $\mathcal O(MQ_S)$
work. The contribution $u_S(\bx_l)$ is then computed for each of the
$N$ targets using $\mathcal O(NQ_S)$ operations.
\end{enumerate}

Since the number of quadrature nodes $Q_S$ is determined by $k$, this
part of the algorithm has a net operation count of $\mathcal
O(nk)$. When $k$ is small, the procedure outlined above is sufficient
for a fast algorithm. If $k$ were large, with $n = \mathcal O(k)$, a
more sophisticated scheme would be required, using either fast
multipole or \emph{butterfly} ideas
\cite{greengard_huang_etc,wideband2d,oneil,rokhlin}.  This would take
us a bit far afield, so we restrict our attention to problems where $n
\gg k$. In that regime, the preceding analysis yields a fast $\mathcal
O (n \log n)$ algorithm for the evaluation of Helmholtz potentials and
fields (for fixed values of $k$ and $\alpha$) with impedance
half-space boundary conditions.

Tables~\ref{tab-timings1} and~\ref{tab-timings2} contain timing
results for randomly located sources and targets that are either
well separated from the interface $y=0$ or located near the
interface. The parameters $N$, $M$, and $M_{image}$ are as
previously defined. The data in the columns labeled $t_{FMM}$,
$t_{SOM}$, and $t_{tot}$ represent the time spent evaluating the
free-space FMM, the Sommerfeld-like integrals, and the total time.
All experiments were performed using Fortran 90 on a laptop with a
$2.53$ GHz Intel Core 2 Duo and 8GB of RAM. The code is not highly
optimized, but the timings do demonstrate the approximately linear
scaling of our algorithm. Straightforward speedups are possible
(e.g. optimizing the location and number of images required,
changing the quadrature used in the Sommerfeld-like integral,
etc.). The impedance Green's functions are calculated to absolute
precision $\epsilon = 10^{-10}$. Note that more images are required
when the sources are located near to the interface $y=0$, i.e. the
data for $M_{image}$ are larger in Table~\ref{tab-timings2} than in
Table~\ref{tab-timings1}. Also, the time required for the
Sommerfeld integral portion of the calculation is greater in
Table~\ref{tab-timings2} than in~\ref{tab-timings1}. This is
because when all the sources and targets are nearer to the
interface $y=0$, fewer quadrature nodes can be used in the
Sommerfeld integral because the integrand $\varphi_{k,\alpha,C}$ is
less oscillatory when $\lambda < k$.

For all the sources in Table~\ref{tab-timings1}, the images could
be omitted as discussed in Remark~\ref{rem-c}. The same holds true
for some portion of the sources in Table~\ref{tab-timings2}.

\subsection{Quadratures for layer potentials}
\label{sec-qbx}

In this section we briefly describe the quadrature scheme we use to
compute singular and weakly singular layer potentials of the type
that appear in (\ref{eq-slprep}) and (\ref{eq-neu}).
For concreteness, we consider the single layer potential $S:
L_2(\Gamma) \to L_2(\Gamma)$:
\begin{equation}
\label{eq-slp}
S[\sigma](\bx) = \int_\Gamma  g_{k,\alpha}(\bx,\by) \, \sigma(\by) \, d\by(s),
\end{equation}
where $\bx \in \Gamma$. Quadrature design for layer potentials is a
well-studied problem and there exist many methods for the numerical
evaluation such integrals. When the discretization is based upon
using equispaced points with respect to some underlying
parameterization, corrected trapezoidal rules are very effective
\cite{alpert,kapur,martinsson}. In the two-dimensional case, these
schemes can integrate functions with logarithmic singularities to
high precision. For each target point $\bx_j \in \Gamma$, they
simply add corrections to the quadrature weights for points in the
vicinity of $\bx_j$. More precisely, by adding $2k$ corrections, it
is possible to obtain errors of the order $O(h^k)$, where $h$ is
the underlying mesh spacing. Better suited for adaptivity are
methods based on piecewise high-order approximations of the density
$\sigma$ \cite{bremer,kolm,martinsson,helsing}. In this work, we
have chosen to use the recently developed QBX method (Quadrature By
Expansion) \cite{klockner}, which requires only a high-order
quadrature rule for {\em smooth} functions. While the method, its
analysis, and its implementation are somewhat intricate, the idea
is simple.

\begin{figure}[t]
  \centering
\begin{subfigure}{.49\textwidth}
    \centering
  \begin{tabular}{|rr|rrrr|}\hline
  $N$ & $M$ & $M_{image}$ & $t_{FMM}$ & $t_{SOM}$ & $t_{tot}$ \\ \hline
  100   & 100   &  12,900  &  0.56 & 0.06 &  0.68 \\
  200   & 200   &  25,800  &  1.12 & 0.12 &  1.38 \\
  400   & 400   &  51,600  &  2.29 & 0.25 &  2.80 \\
  800   & 800   &  103,200 &  4.78 & 0.49 &  5.79 \\
  1,600  & 1,600  &  206,400 &  9.68 & 0.98 & 11.70 \\
  3,200  & 3,200  &  412,800 & 19.43 & 1.97 & 23.46 \\
  6,400  & 6,400  &  825,600 & 39.21 & 3.93 & 47.27 \\
\hline
  \end{tabular}
  \caption{Sources and targets in $(-1,1) \times (2,3)$.}
\label{tab-timings1}
\end{subfigure} 
\begin{subfigure}{.49\textwidth}
    \centering
  \begin{tabular}{|rr|rrrr|}\hline
  $N$ & $M$ & $M_{image}$ & $t_{FMM}$ & $t_{SOM}$ & $t_{tot}$ \\ \hline
  100   & 100   & 17,268   & 0.71  & 0.03 &  0.77 \\
  200   & 200   & 34,920   & 1.45  & 0.05 &  1.57 \\
  400   & 400   & 70,352   & 2.95  & 0.10 &  3.18 \\
  800   & 800   & 140,720  & 6.08  & 0.21 &  6.53 \\
  1,600  & 1,600  & 281,696  & 12.32 & 0.41 & 13.22 \\
  3,200  & 3,200  & 562,176  & 25.19 & 0.83 & 26.98 \\
  6,400  & 6,400  & 1,122,960 & 51.22 & 1.65 & 54.79 \\ \hline
  \end{tabular}
  \caption{Sources and targets in $(-1,1) \times (0,1)$.}
  \label{tab-timings2}
\end{subfigure}
\caption{Timing results for $N$ targets and $M$ sources uniformly
    randomly distributed, $k=10.2$, and $\alpha=2.04$.}
\end{figure}

For each target point $\bx_j \in \Gamma$ with normal $\bn_j$, let
$\bc_j$ denote a nearby off-surface point, say at
\[
    \bc_j = \bx_j + r \bn_j .
\]
Assuming there are no sources in the disk of radius $r$ about $\bc_j$,
$S[\sigma](\bx)$ satisfies the homogeneous Helmholtz equation in that
disk. Using the standard separation of variables solution to the
Helmholtz equation, we may therefore write
\begin{equation} \label{qbxexpand}
S[\sigma](\bx) \approx \sum_{l=-p}^p a_l \, J_l(k\rho) \, e^{-il\theta} ,
\end{equation}
where $(\rho,\theta)$ are the polar coordinates of $\bx$ with respect
to the expansion center $\bc_j$ and $J_l$ denotes the Bessel function
of the first kind of order $l$.  Rather than evaluate (\ref{eq-slp})
directly, in QBX one evaluates the expansion (\ref{qbxexpand})
instead.  It is shown in \cite{klockner,klockner2} that the error in
evaluating this local expansion is of the order $\mathcal
O(r^{p+1}~+~(\frac{h}{r})^{2q})$, where $h$ is the discretization
spacing and $q$ is the order of the underlying quadrature rule for
{\em smooth} functions. The method converges like a $p^\text{th}$-order
accurate scheme until the $( \frac{h}{r} )^{2q}$ error dominates.
This latter term can be made arbitrarily small.

QBX is particularly useful when applied to layer potentials with the
impedance half-space Green's function.  If the boundary $\Gamma$ of
the obstacle $\Omega$ is very close to the $x$-axis, say within
$\epsilon$ of the interface, then the near-field images are located in
the lower-half plane at a distance only $2\epsilon$ away. The
proximity of these images causes difficulties for standard quadrature
techniques because of the near singularities they induce.  The QBX
expansion centers $\bc_j$, however, may be placed in the {\em interior} 
of $\Omega$, and the integrals are computed easily.

\section{Numerical results}
\label{sec_numerical}

We now have the necessary machinery needed to solve non-trivial
scattering problems using an integral equation formulation. We
discretize the (smooth) scatterer at equispaced points with respect to
the underlying parameterization, and, unless otherwise noted, use a
$16^{\text{th}}$-order QBX quadrature scheme to evaluate layer
potentials. We solve the discretized integral equation using the
iterative method GMRES, accelerated by the fast algorithm described in
Section~\ref{sec-fast}.

The numerical examples in this section illustrate the accuracy of the
computation of $g_{k,\alpha}$, as well as the robustness of our hybrid
representation with respect to the location of the
sources and targets (in particular, when both source and target are
located near the interface $y=0$).

\begin{figure}[p]
\centering
\begin{subfigure}[b]{.45\textwidth}
    \centering
    \includegraphics[width=\textwidth]{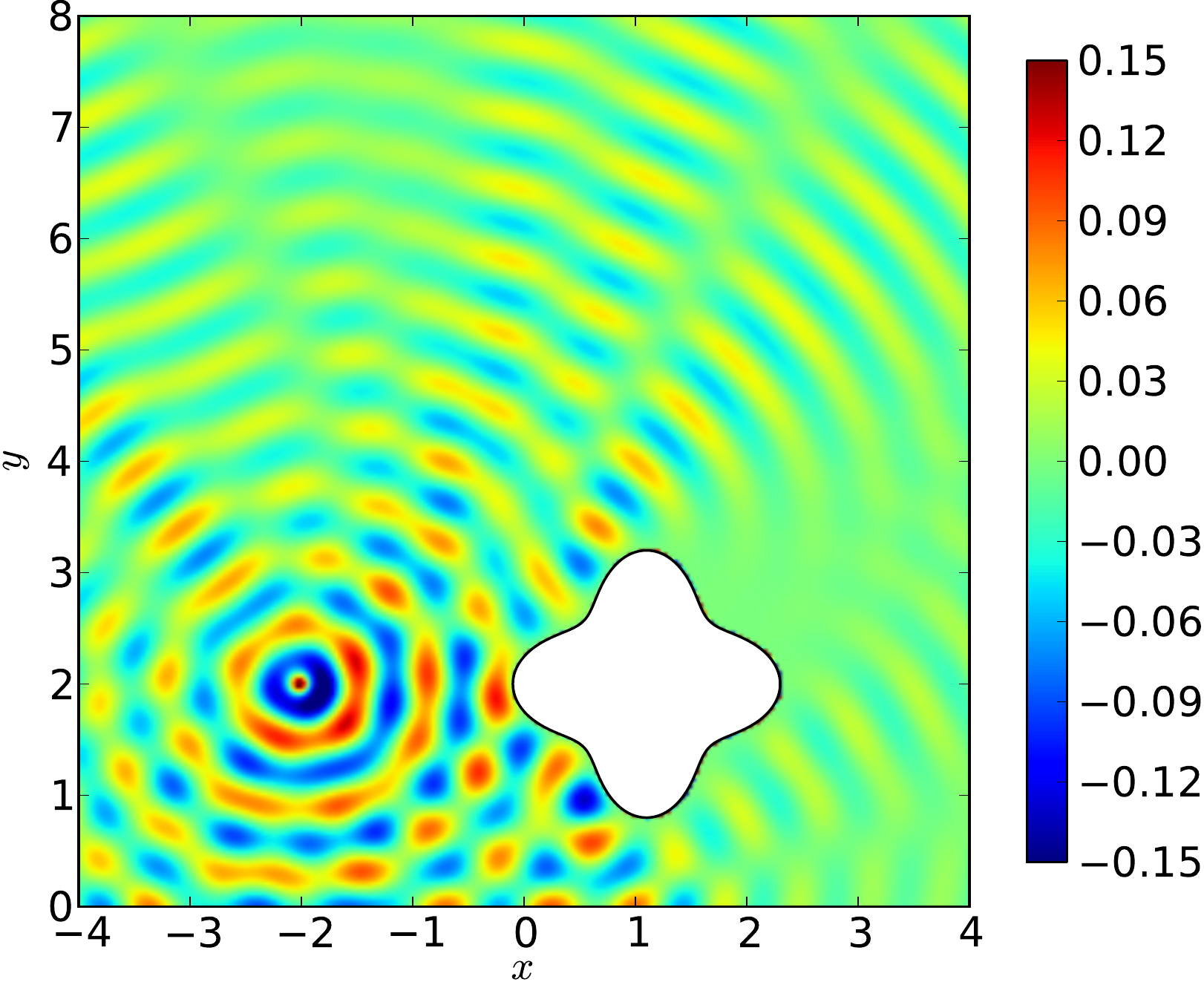}
    \caption{Exterior Dirichlet, $\delta=0.8$.}
    \label{fig-dir1}
\end{subfigure} \qquad \qquad
\begin{subfigure}[b]{.45\textwidth}
    \centering
    \includegraphics[width=\textwidth]{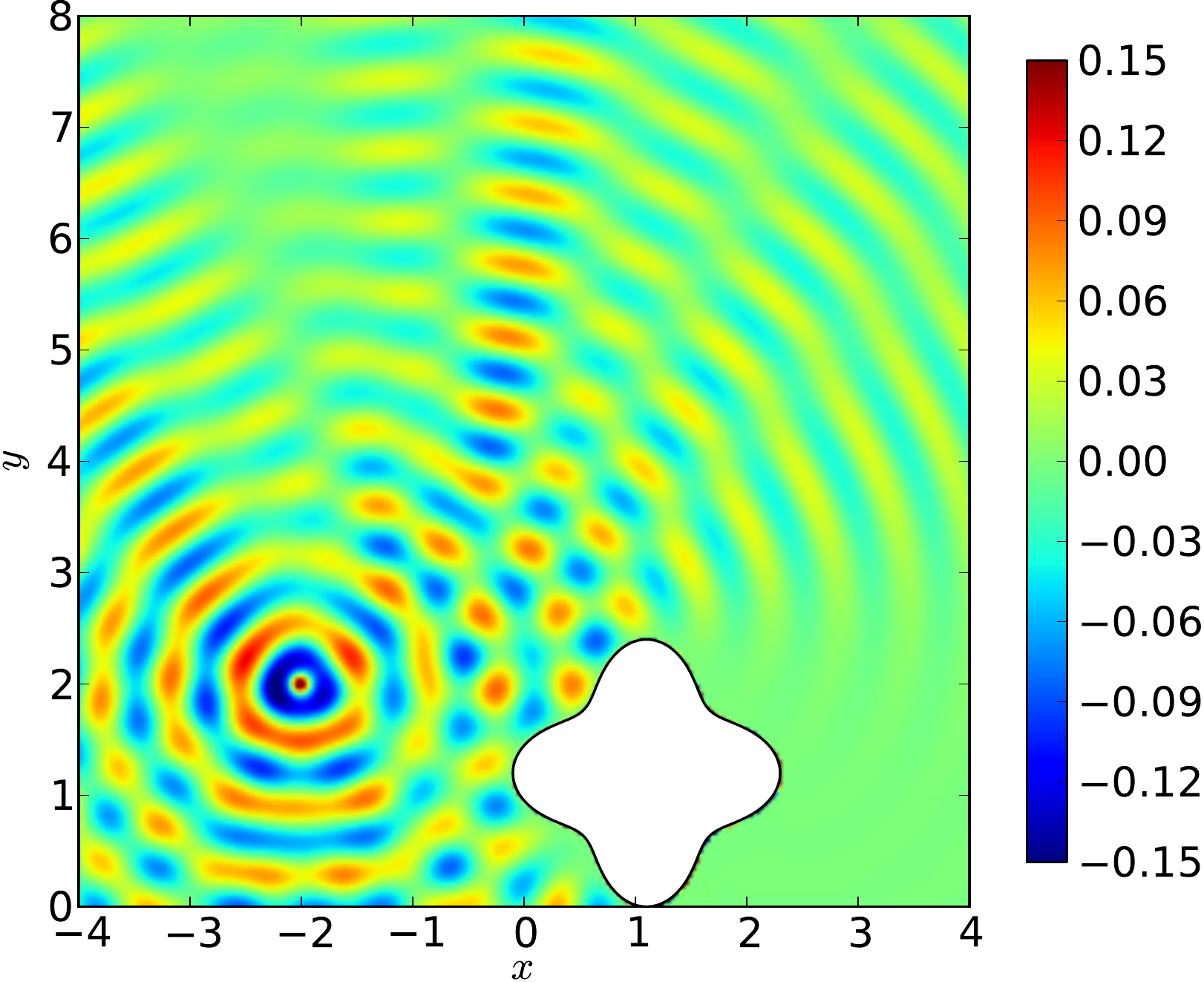}
    \caption{Exterior Dirichlet, $\delta=10^{-3}$.}
    \label{fig-dir2}
\end{subfigure} \\
\vspace{.5in}
\begin{subfigure}[b]{.45\textwidth}
    \centering
    \includegraphics[width=\textwidth]{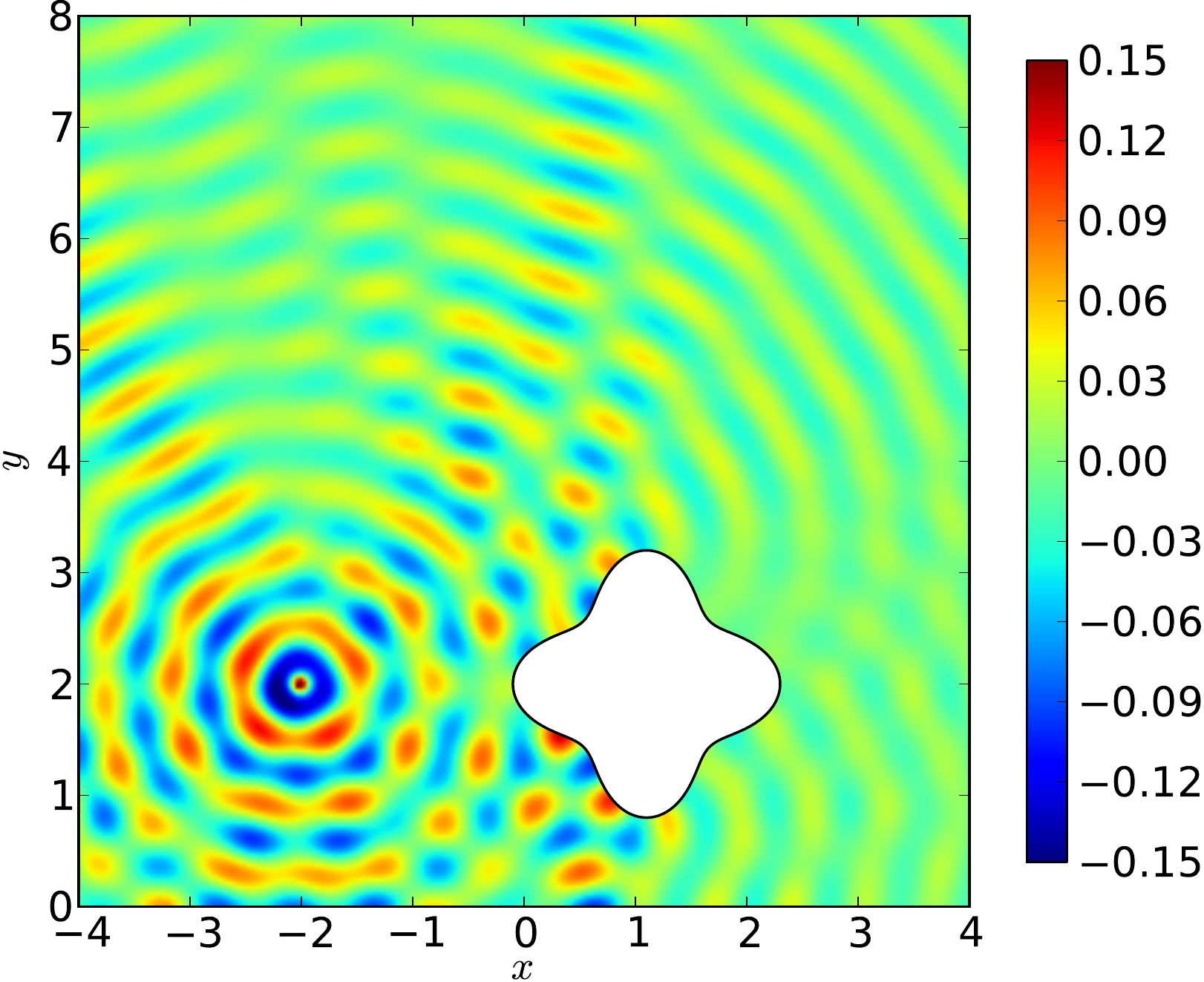}
    \caption{Exterior Neumann, $\delta=0.8$.}
    \label{fig-neu1}
\end{subfigure} \qquad \qquad
\begin{subfigure}[b]{.45\textwidth}
    \centering
    \includegraphics[width=\textwidth]{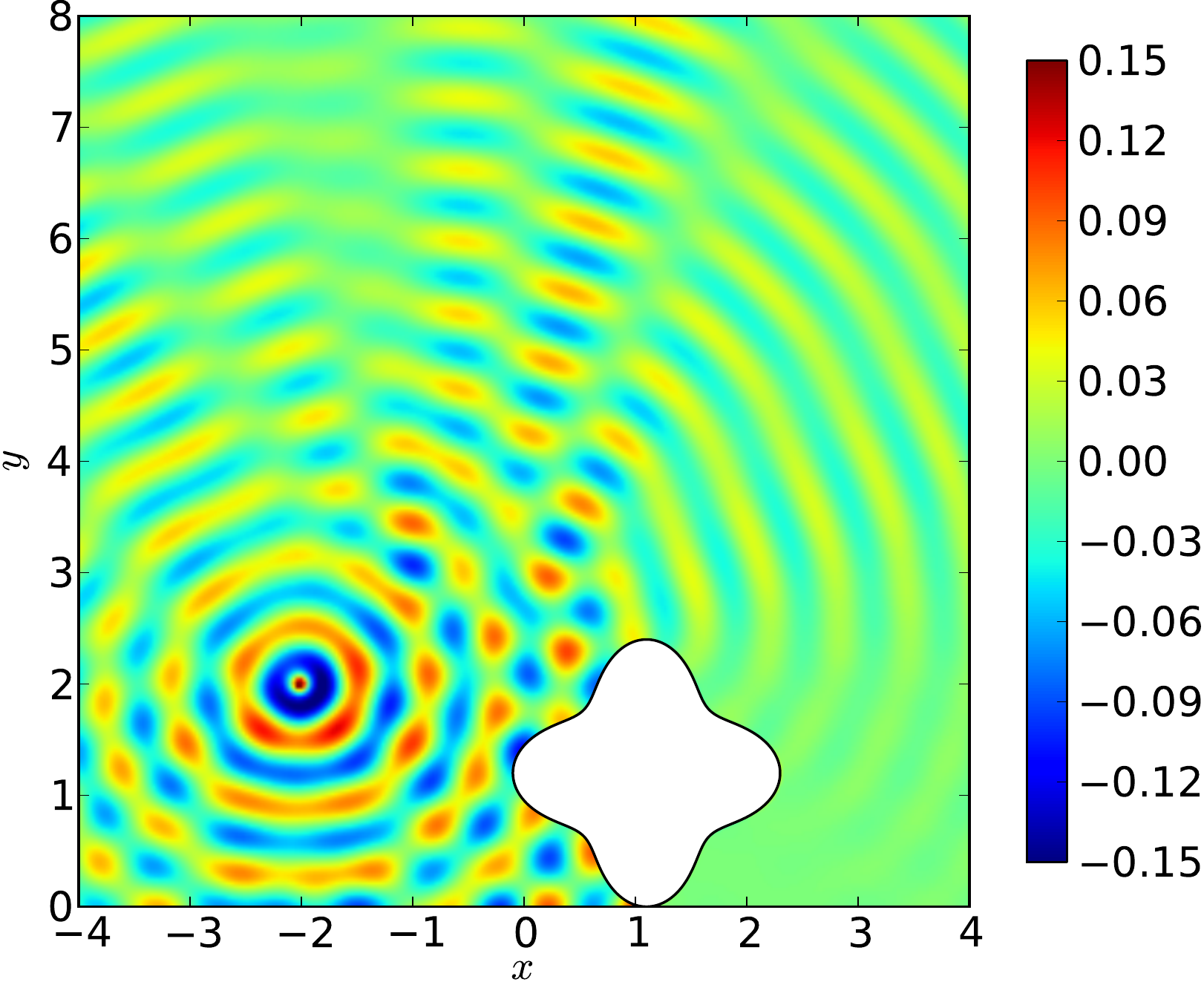}
    \caption{Exterior Neumann, $\delta=10^{-3}$.}
    \label{fig-neu2}
\end{subfigure} \\
\caption{For $k=10.2$ and $\alpha = 2.04$, we show the real part of
  the total field for exterior Dirichlet and Neumann
  problems.  The incoming field in all cases is due to a unit 
  strength impedance point source located at $(-2, 2)$.  The 
  distance of the scatterer from the interface is listed below 
  each figure.}
\label{fig-main}
\end{figure}

We first solve the exterior
(sound-soft) Dirichlet problem with impedance interface conditions:
\begin{equation}
\begin{split}
(\triangle + k^2) u^{tot}(\bx) = \delta(\bx-\bx_0) & \qquad \text{in } P \setminus \Omega ,\\
\frac{\partial u^{tot}}{\partial n} -  i \alpha u^{tot} = 0 & \qquad 
    \text{on } y=0 ,\\
u^{tot} = 0 & \qquad \text{on } \Gamma,
\end{split}
\end{equation}
where $\Omega$ is an inclusion in the interior of the upper
half-space.  We assume a representation of the scattered field in the
form of a double layer potential:
\begin{equation}
u(\bx) = \int_\Gamma 
\left[ \frac{\partial}{\partial n_y} g_{k,\alpha}(\bx,\by) \right]
 \sigma(\by) \, d\by(s), 
\end{equation}
where as before $u^{tot} = u + u^{in}$.
This representation leads to the integral equation
\begin{equation}
\label{eq-dir}
\frac{1}{2} \sigma(\bx) + \int_\Gamma
\left[ \frac{\partial}{\partial n_y} g_{k,\alpha}(\bx,\by) \right]
 \sigma(\by) \, d\by(s) = -u^{in}(\bx),
\end{equation}
which is a Fredholm equation of the second-kind, with the integral
interpreted in the principal value sense.  The equation
(\ref{eq-dir}), like (\ref{eq-neu}), is invertible except for a
countable sequence of spurious resonances.  Rather than use a more
complicated, resonance-free representation, we assume $k$ is not a
spurious resonance.  See \cite{chandler-wilde} for a thorough
discussion of integral equations in the impedance scattering context.
Figures~\ref{fig-dir1} and~\ref{fig-dir2} show the total potential 
in the exterior of the scatterer $\Omega$ with zero Dirichlet 
conditions. In
these, and the subsequent Neumann scattering examples, the boundary of
the scatterer, $\Gamma$, is described by the curve $\gamma : [0,2\pi]
\to \mathbb R^2$, where $\gamma(t) = (x(t), y(t))$ and
\begin{equation}
\label{eq-curve}
\begin{split}
x(t)&=1.1  + (1 + 0.2 \cos{4t} ) \cos{t} ,\\
y(t)&=1.2 + \delta + (1 + 0.2 \cos{4t} ) \sin{t}.
\end{split}
\end{equation}
The magnitude $\delta$ of the distance to the interface $y=0$ is
indicated below below each figure. Accuracy results are given in
Table~\ref{tab-dirneu}. Note that more discretization points were
required when $\delta=10^{-3}$, i.e. when the scatterer is close to
the impedance boundary. An adaptive discretization scheme, instead of
a global trapezoidal rule, would result in many fewer nodes. A
slightly higher-order quadrature, $20^{\text{th}}$-order,
was used in the $\delta=10^{-3}$ case.

\begin{table}[t]
  \centering
  \begin{tabular}{|l|cccc|}\hline
& Dirichlet, $\delta=0.8$ & Dirichlet, $\delta=10^{-3}$
& Neumann, $\delta=0.8$ & Neumann, $\delta=10^{-3}$ \\ \hline
Number of points              & 500     & 1500     & 500     & 1500 \\
Rel. error at target $u$      & 0.10E-9 & 0.15E-10 & 0.22E-9 & 0.14E-11 \\
Rel. $L_2$ error of $\sigma$  & 0.29E-9 & 0.32E-10 & 0.40E-9 & 0.98E-8 \\
$L_2$ norm of $\sigma$        & 2.13    & 2.95     & 4.03    & 3.50 \\
\hline
  \end{tabular}
  \caption{Accuracy results.} 
\label{tab-dirneu}
\end{table}

In the second set of experiments, we solve the exterior
(sound-hard) Neumann problem with impedance half-space conditions:
\begin{equation}
\begin{split}
(\triangle + k^2) u^{tot}(\bx) = \delta(\bx-\bx_0) & \qquad \text{in } P 
\setminus \Omega ,\\
  \frac{\partial u^{tot}}{\partial n} - i \alpha u^{tot} = 0 &
  \qquad \text{on } y=0 ,\\ \frac{\partial u^{tot}}{\partial n} = 0 &
  \qquad \text{on } \Gamma.
\end{split}
\end{equation}
As discussed in the introduction, we represent the scattered field by a
single layer potential (\ref{eq-slprep}),
yielding the integral equation (\ref{eq-neu}).
Figures~\ref{fig-neu1} and~\ref{fig-neu2}
show the total potential in the exterior of the scatterer $\Omega$.

It is worth repeating that the interface $y=0$ does \emph{not} need to
be discretized because of the use of the impedance Green's function
$g_{k,\alpha}$.  The number of discretization points and quadrature
parameters for our examples were not carefully chosen. They were
simply chosen to be sufficiently fine to yield an estimated precision
in the scattered field of $10^{-10}$. The impedance Green's function was
evaluated to absolute precision $10^{-11}$ and the 
GMRES residual tolerance was set to $10^{-10}$. Our
goal here is simply to demonstrate the accuracy and robustness of our
hybrid representation of the Green's function.

The accuracy of the solvers is tested by two different methods. First,
a unit strength source is placed at a point $\bx_I$ in the
\emph{interior} of the scatterer $\Omega$.  The function
$g_{k,\alpha}(\bx,\bx_I)$ is then a known exact solution satisfying the
homogeneous Helmholtz equation in $P\setminus \Omega$.  Using the
Dirichlet or Neumann data corresponding to $g_{k,\alpha}(\bx,\bx_I)$,
we can solve the boundary integral equations (\ref{eq-neu}) and
(\ref{eq-dir}) and compare the computed solution at target points in
$P\setminus \Omega$ with the known exact solution.  For the true
scattering problem, the source is exterior to $\Omega$ and an exact
solution is not available in closed form.  In this case, we carry out
a self-consistent convergence test. That is, the number of
discretization points is doubled and the relative $L_2$ error is
calculated with the finer grid solution used as the reference
solution.  The error in the potential at a target point in the
exterior of $\Omega$ is calculated in the same fashion.  Results of
the self-consistent convergence tests are shown in
Table~\ref{tab-dirneu}.  The location of the target point in the
exterior at which the value of the potential was computed is $(0,5)$
in all examples.

\subsection{Scattering from a locally perturbed half-space}
\label{perturb}

\begin{figure}[t]
\centering
\begin{subfigure}[b]{.6\textwidth}
    \centering
    \includegraphics[width=\textwidth]{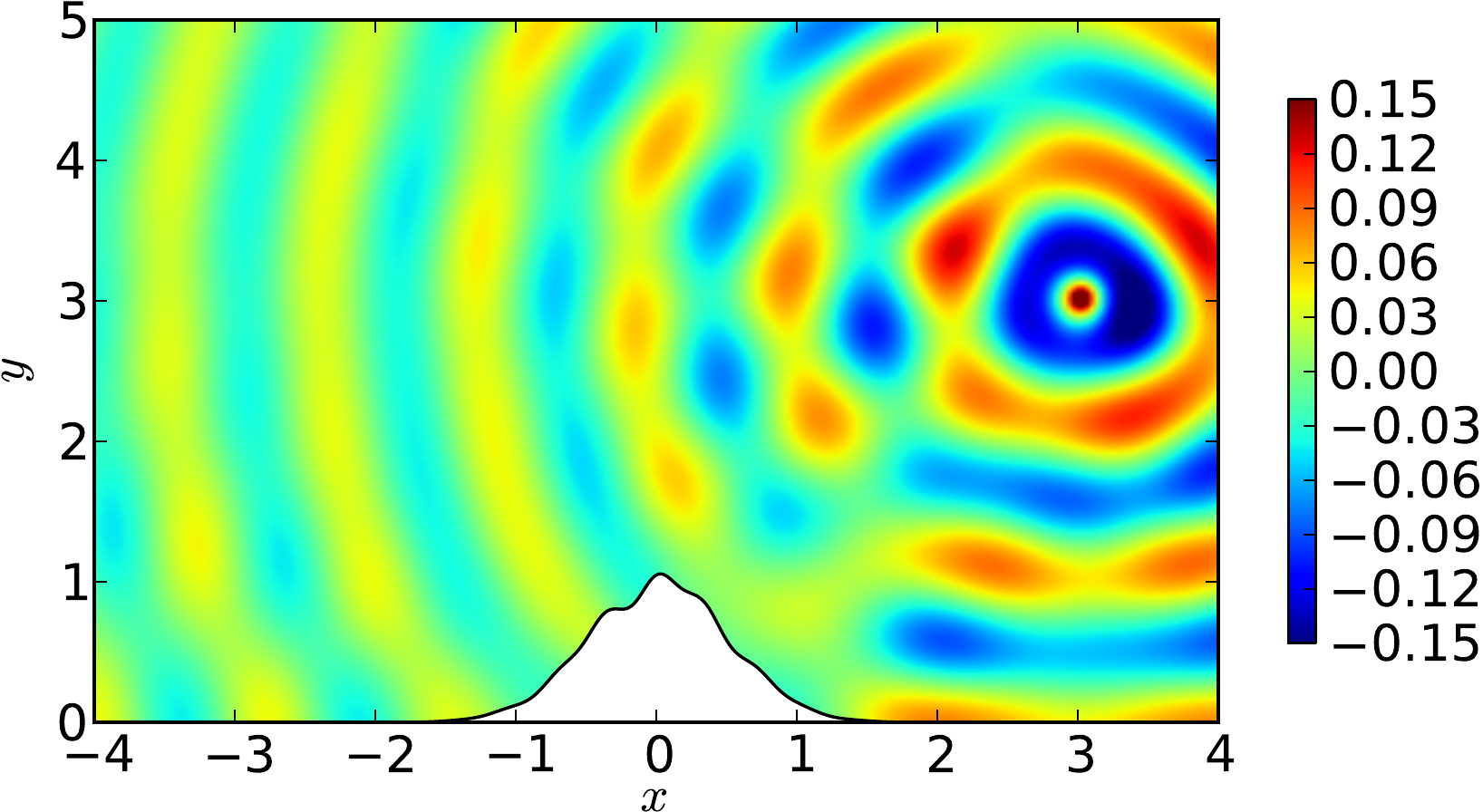}
    \caption{Impedance scattering from a half-space perturbation.}
    \label{subfig-surf}
\end{subfigure} \qquad \qquad
\begin{subfigure}[b]{.3\textwidth}
    \centering
    \begin{tabular}{|l|c|}\hline
    Number of points                & 4,000 \\ 
    Rel. $L_\infty$ error of $u$   &  0.40E-9 \\
    Rel. $L_2$ error of $\sigma$    & 0.51E-8  \\
    $L_2$ norm of $\sigma$    & 0.76 \\ \hline
    \end{tabular}
    \vspace{.75in}
    \caption{Accuracy results.}
    \label{tab-surf}
\end{subfigure}
\caption{The real part of 
the total potential is plotted for $k=5.7$ and $\alpha = .855$ over a
perturbed impedance half-space.
The incoming field is due to a unit strength point source
located at $(3, 3)$. Accuracy results are obtained by doubling the 
number of discretization points, for a test target point located at $(-2,4)$.
\label{fig-surf}}
\end{figure}

Lastly, we demonstrate the use of the impedance Green's function in
solving the problem of scattering from a half-space which contains a
local perturbation. A thorough analysis can be found
in~\cite{chandler-wilde-peplow}. Here, we simply consider one such
example.  We let $\Omega$ denote the region below the curve~$\gamma_1:
(-\infty,\infty) \to \mathbb R^ 2$, where $\gamma_1(t) = (x(t), y(t))$
and
\begin{equation}
\begin{split}
x(t) &= t, \\
y(t) &= \left( 1+0.05 \left( \sin 8.79 t + \cos 16.96 t+\sin1.88 t
\right) \right)
\, e^{-2 t^2}.
\end{split}
\end{equation}
For values of $t$ outside the interval $[-4,4]$, $|y(t)| < 2.0 \times
10^{-14}$.  Therefore, $\gamma_1$ is only discretized 
for $t \in [-4,4]$, outside of which
it is indistinguishable (to machine precision) from the interface
$y=0$. Thus, outside of $t \in [-4,4]$, the impedance condition is
automatically satisfied to machine precision
if the appropriate Green's function is used to represent the solution.
On the half-space perturbation $\Gamma = \{\gamma_1(t), t \in [-4,4]\}$,
we must enforce the desired impedance condition:
\begin{equation}
\frac{\partial u^{tot}}{\partial n} - i \alpha u^{tot} = 0 \qquad
\text{on } \Gamma,
\end{equation}
where $n=(y',-x')$ is a unit normal vector pointing into the lower
half space $\Omega$.  In short, we seek to solve the boundary value
problem:
\begin{equation}
\begin{split}
(\triangle + k^2) u^{tot}(\bx) = \delta(\bx-\bx_0) & \qquad \text{in } P\setminus\Omega, \\
  \frac{\partial u^{tot}}{\partial n} - i \alpha u^{tot} = 0 &
  \qquad \text{on } \Gamma.
\end{split}
\end{equation}

We assume an integral representation for the scattered field $u$:
\begin{equation}
u(\bx) = \int_\Gamma g_{k,\alpha}(\bx,\by) \,
 \sigma(\by) \, d\by(s), 
\end{equation}
which results in the integral equation 
\begin{equation}
\label{eq-inteq-surf}
\frac{1}{2} \sigma(\bx) + \int_\Gamma
\left[ \frac{\partial}{\partial n_x} g_{k,\alpha}(\bx,\by) \right]
 \sigma(\by) \, d\by(s) - i \alpha 
\int_\Gamma g_{k,\alpha}(\bx,\by) \, \sigma(\by) \, d\by(s)
 = -\frac{\partial u^{in}(\bx)}{\partial n_x} + i \alpha u^{in}(\bx).
\end{equation}
This is a simplification of the integral equation approach
discussed in \cite{chandler-wilde-peplow}, which addresses the more
general case of a half-space perturbation that is allowed to
deviate both above and below the interface $y=0$. For simplicity,
we have assumed that the perturbation $\gamma_1$ satisfies
$y(t)\geq 0$.

Even though the curve $\gamma_1$ in this half-space scattering
example comes arbitrarily close to the interface $y=0$, we obtain
high accuracy (using a $16^\text{th}$ order QBX scheme with the
trapezoidal rule as the underlying smooth rule). The reason for
this is that the density $\sigma$ decays exponentially as $|t|$
increases and is as smooth as the perturbation itself. Unlike the
closely-touching Dirichlet and Neumann scattering problems, the
half-space perturbation problem is not physically ill-conditioned.

Figure~\ref{fig-surf} shows the real part of the total potential
$u^{tot}$. The incoming potential is due to an impedance point
source located at $(3,3)$, and the error in the total potential is
calculated at $(-2,4)$. Errors are calculated in the same manner as
for the Dirichlet and Neumann scattering examples, using a
self-consistent convergence test.

\begin{figure}[t]
\centering
\begin{subfigure}[b]{.48\textwidth}
    \centering
    \includegraphics[width=\textwidth]{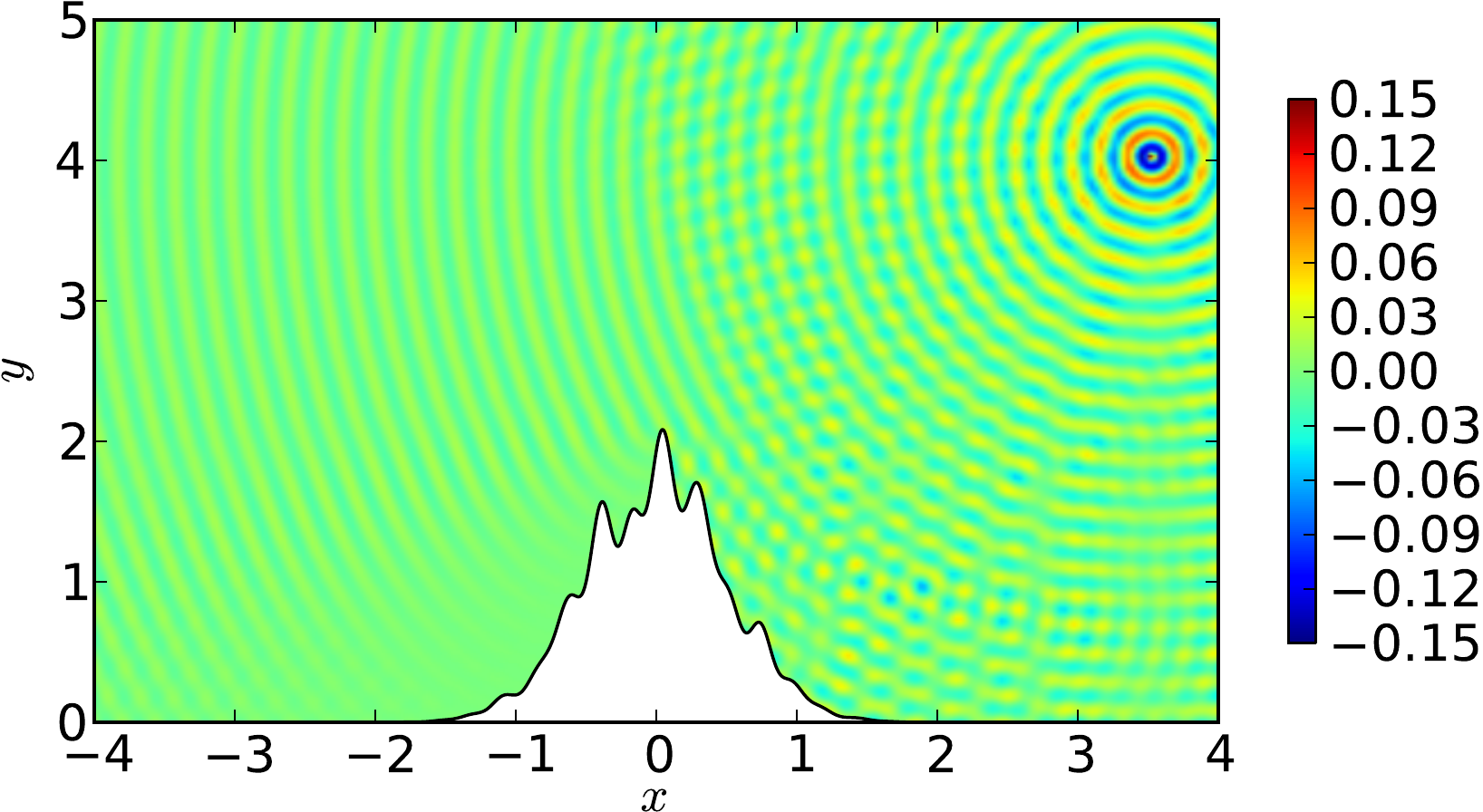}
    \caption{Impedance scattering from a half-space perturbation.}
    \label{subfig-surf3a}
\end{subfigure} 
\begin{subfigure}[b]{.48\textwidth}
    \centering
    \includegraphics[width=\textwidth]{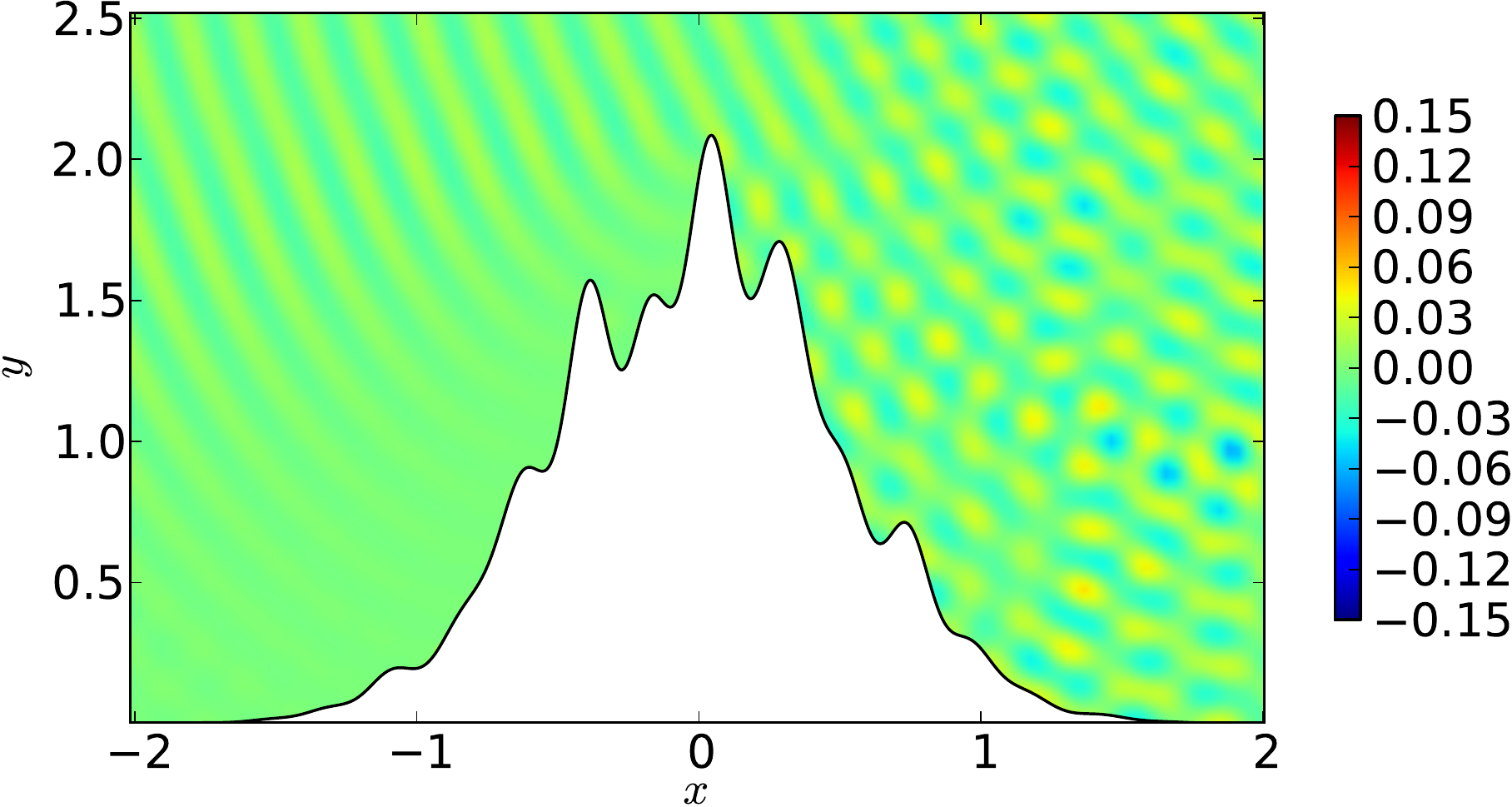}
    \caption{A closer view.}
    \label{subfig-surf3b}
\end{subfigure}
\caption{The real part of 
the total potential is plotted for $k=31.7$ and $\alpha = 5.389$ over a
perturbed impedance half-space.
The incoming field is due to a unit strength point source
located at $(3.5, 4)$. The perturbation was discretized with $20,000$
points, and a relative $L_2$ error of $10^{-5}$ was obtained
in the density $\sigma$ and $10^{-7}$ in $u^{tot}$
at $(-2,5)$. The deviation is approximately $60$ wavelengths long.
\label{fig-surf2}}
\end{figure}

Figure~\ref{fig-surf2} shows the real part of the total potential
$u^{tot}$ for a higher value of the wavenumber $k$ than
Figure~\ref{fig-surf}, and a more complicated perturbation
$\gamma_2(t) = (x(t), y(t))$ with
\begin{equation}
\begin{split}
x(t) &= t, \\
y(t) &= \left( 1.75+\frac{1}{6} \left( \sin 8.79 t + \cos 16.96
t+\sin 27.02 t + \frac{1}{7} \cos 32.67 t \right) \right) \, e^{-2
t^2}.
\end{split}
\end{equation}
The incoming potential is due to an impedance point source slightly
farther away, located at $(3.5,4)$, and the error in the total
potential is calculated at $(-2,5)$. Errors are calculated as
previously discussed.

\section{Conclusions}
\label{sec_conclusions}

We have derived a new formula for the half-space Helmholtz Green's
function satisfying impedance boundary conditions in two
dimensions. The representation (\ref{eq:hybrid_rep}) consists of a
free-space Helmholtz Green's function in the upper half-space, a
short segment of images in the lower half-space with {\em real}
coordinates, and a rapidly converging Sommerfeld-like integral.
Unlike the method of complex images, it is straightforward to
accelerate with an FMM using a modest number of discrete image
charges. The impedance Green's function is easily evaluated to full
double-precision accuracy, insensitive (within logarithmic factors)
to the location of the source and target, using a modest number of
operations that depends only on $k$ and $\alpha$.

Although we have carried out a detailed derivation only for the
two-dimensional case, the results are nearly identical in three
dimensions. Using the three-dimensional Green's function
$g_k(\bx,\by) = -\frac{e^{ik|\bx-\by|}}{4\pi |\bx-\by|}$,
it is straightforward to show that 
\begin{equation}
\begin{split}
g_{k,\alpha}(\bx,\bx_0) =& \ g_k(\bx,\bx_0) +
g_k(\bx,\bx_0-2 y_0 \hat\by) +2i\alpha
\int_0^{C} g_k(\bx,\bx_0-(2 y_0 + \eta)\hat\by) \, 
 e^{-i\alpha\eta} \, d\eta \\
 &+\frac{i\alpha}{2\pi} \int_{-\infty}^\infty 
\frac{e^{-\sqrt{\lambda^2-k^2}(y+y_0)}} {\sqrt{\lambda^2-k^2} }
\frac{e^{-(\sqrt{\lambda^2-k^2}+i\alpha)C}}
{\sqrt{\lambda^2-k^2} +i\alpha}
J_0(\lambda r) \, d\lambda,
\end{split} \label{eq:hybrid_rep2d}
\end{equation}
where $J_0$ is the zeroth-order Bessel function of the first
kind.

A more complicated extension of our approach is to the evaluation
of the layered medium Green's function, where each layer has a
distinct Helmholtz coefficient and various (application-dependent)
continuity conditions are imposed across each layer. In this case,
an explicit real image structure is not available via the Laplace
transform. We are currently investigating a semi-numerical approach
that appears promising. We are also developing hybrid
representations for the full Maxwell equations with impedance
boundary conditions \cite{oneil-sifuentes}.

The reader may have noted that formulas
(\ref{eq_im2}) and (\ref{eq:hybrid_rep}) make use of  a discrete
image at $(x_0,-y_0)$ with the same sign as the original source. 
As $\alpha \rightarrow 0$, the impedance
condition becomes a sound-hard condition and the discrete image
is all that remains, as one would expect. 
As $|\alpha| \rightarrow \infty$, however, the limit of the 
impedance condition is a Dirichlet condition, while
the image formula diverges. One should expect to see a simple image
source at $(x_0,-y_0)$ with the opposite sign. 
This is easy to accomplish by 
using a complementary image formula.  
Rather than (\ref{eq_mu}), one can write
\begin{equation}\label{eq_eta2}
\begin{split}
\int_0^\infty e^{-\sqrt{\lambda^2-k^2}\eta} \, \tau(\eta) \, d\eta &=
\frac{\sqrt{\lambda^2-k^2} - i\alpha}{\sqrt{\lambda^2-k^2} + i\alpha} \\
&=-1 + \frac{2 \sqrt{\lambda^2-k^2} }{\sqrt{\lambda^2-k^2} + i\alpha} \, .
\end{split} 
\end{equation}
This has the desired asymptotic behavior as $\alpha \rightarrow \infty$.
In applications $\alpha$ is typically $\mathcal O(|k|)$, and the 
formula presented in Section~\ref{sec:hybrid} serves its purpose well.


\end{document}